\font\cmrfootnote=cmr10 scaled 750
\font\cmrhalf=cmr10 scaled \magstephalf

\font\cmrscfootnote=cmr7 scaled 750
\font\cmrschalf=cmr7 scaled \magstephalf

\font\cmrscscfootnote=cmr5 scaled 750
\font\cmrscschalf=cmr5 scaled \magstephalf

\font\mitfootnote=cmmi10 scaled 750
\font\mithalf=cmmi10 scaled \magstephalf

\font\mitscfootnote=cmmi7 scaled 750
\font\mitschalf=cmmi7 scaled \magstephalf

\font\mitscscfootnote=cmmi5 scaled 750
\font\mitscschalf=cmmi5 scaled \magstephalf

\font\cmsyfootnote=cmsy10 scaled 750
\font\cmsyhalf=cmsy10 scaled \magstephalf

\font\cmsyscfootnote=cmsy7 scaled 750
\font\cmsyschalf=cmsy7 scaled \magstephalf

\font\cmsyscscfootnote=cmsy5 scaled 750
\font\cmsyscschalf=cmsy5 scaled \magstephalf

\font\cmexfootnote=cmex10 scaled 750
\font\cmexhalf=cmex10 scaled \magstephalf

\font\cmexscfootnote=cmex10 scaled 750
\font\cmexschalf=cmex10 scaled \magstephalf

\font\cmexscscfootnote=cmex7 scaled 750
\font\cmexscschalf=cmex7

\def\mathfootnote{\textfont0=\cmrfootnote \textfont1=\mitfootnote \textfont2=\cmsyfootnote \textfont3=\cmexfootnote
        \scriptfont0=\cmrscfootnote \scriptscriptfont0=\cmrscscfootnote \scriptfont1=\mitscfootnote \scriptscriptfont1=\mitscscfootnote
        \scriptfont2=\cmsyscfootnote \scriptscriptfont2=\cmsyscscfootnote \scriptfont3=\cmexscfootnote \scriptscriptfont3=\cmexscscfootnote}
\def\mathhalf{\textfont0=\cmrhalf \textfont1=\mithalf \textfont2=\cmsyhalf \textfont3=\cmexhalf
        \scriptfont0=\cmrschalf \scriptscriptfont0=\cmrscschalf \scriptfont1=\mitschalf \scriptscriptfont1=\mitscschalf
        \scriptfont2=\cmsyschalf \scriptscriptfont2=\cmsyscschalf \scriptfont3=\cmexschalf \scriptscriptfont3=\cmexscschalf}

\let\mf=\mathfootnote

\font\Bbb=msbm10 
\font\Bbbfootnote=msbm7 scaled \magstephalf

\def\outlin#1{\hbox{\Bbb #1}}
\def\outlinfootnote#1{\hbox{\Bbbfootnote #1}}

\font\call=cmsy10 

\def\cal{\call}
\font\calfoot=cmsy7

\font\small=cmr5                       
\font\notsosmall=cmr7
\font\cmreight=cmr8
 
\font\bfeight=cmbx8

\font\cmrhalf=cmr10 scaled \magstephalf
       
\font\two=cmr10 scaled \magstep2

\font\sansfoot=cmss8
\font\sans=cmss10

\font\twosans=cmss10 scaled \magstep2

\font\caps=cmcsc10

\let\mf=\mfootnote

\def\q{\quad}   \def\qq{\qquad}

\def\cut{\hfill\break}  \def\newline{\cut}
\def\h#1{\hbox{#1}}

\def\eqa{\eqalign} \def\eqano{\eqalignno}

 \def\eps{\epsilon}
\def\be{\beta}
  
 \def\L{\Lambda} 
\def\vp{\varphi}
  
\def\o{\omega}   \def\O{\Omega}

\def\calB{{\h{\cal B}}} \def\calH{{\h{\cal H}}}
\let\H=\calH 
\def\calM{{\h{\cal M}}} 
  
 \def\calD{{\h{\cal D}}}
 
\def\calL{{\h{\cal L}}} \def\L{{\calL}} 
 \def\calA{{\h{\cal A}}}

\def\calF{{\h{\cal F}}}

\def\calfootH{{\h{\calfoot H}}}  \def\Hfoot{\calfootH}
\def\calfootC{{\h{\calfoot C}}}

\def\CC{{\outlin C}} \def\RR{{\outlin R}} \def\ZZ{{\outlin Z}} 
  \def\PP{{\outlin P}}

\def\CCfoot{{\outlinfootnote C}}

\chardef\dotlessi="10  
\chardef\inodot="10

\def\polishL{\leavemode\setbox0=\hbox{L}\hbox to\wd0{\hss\char'40L}}

\def\Dom{\h{\rm Dom}}

\def\Aut{\h{\rm Aut}}  \def\aut{\h{\rm aut}}

 \def\Dom{\hbox{\rm Dom}}  \def\Ric{\h{\rm Ric}}

\def\Iso{\h{\rm Iso}}

\def\V{\frac1V}

\def\gHerm{g_{\h{\small Herm}}}

\def\gij{g_{i\bar j}}

\def\oFSc{\o_{\hbox{\small FS},c}}

\def\GL2nR{\h{$GL(2n,\RR)$}}
\def\Sp2nR{\h{$Sp(2n,\RR)$}}

  \def\dz{d{\bf z}} 
\def\dz#1{dz^{#1}}
\def\dzb#1{d\overline{z^{#1}}}  
\def\dbz{d\bar z}

\def\dzidzjb{dz^i\w\dbz^j}

\def\ddt{\frac d{dt}}

\def\part#1{\frac{\partial#1}{\partial t}}

\def\ov{\over}

\def\frac#1#2{{{#1}\over{#2}}}

\def\all{\forall}

\def\sm{\setminus}

\def\precpt
{\hbox{$\mskip3mu\mathhalf\subset\raise0.92pt\hbox{$\mskip-10mu\!\!\!\!
\mathfootnote\subset$}\mskip5mu$}}

\def\supsetnoteq{\hbox{$\mskip3mu\supset\raise-5.97pt
\hbox{$\mskip-10mu\!\!\!\scriptstyle\not=$}\mskip8mu$}}

\def\subsetnoteq{\hbox{$\mskip3mu\subset\raise-5.97pt
\hbox{$\mskip-11mu\!\!\scriptstyle\not=$}\mskip8mu$}}

\def\sseq{\subseteq}

\def\isom{\cong}

\def\*{\star}

\def\w{\wedge}

\def\ii#1{{#1}^{-1}}

\def\D{\Delta}
 
\def\dbar{\bar\partial}
\def\del{\partial} 
\def\ddbar{\partial\dbar}

\def\intM{\int_M}  \def\intm{\int_M}
\def\V{\frac1V}

 \def\ha{\frac12}

\def\i{\sqrt{-1}}

\def\ra{\rightarrow}

\def\arrow#1{\hbox to #1pt{\rightarrowfill}}

\def\thhnotsosmall#1{${\hbox{#1}}^{\hbox{\small th}}$}

\def\MA{Monge-Amp\`ere }

\def\Kahler{K\"ahler }
\def\K{K\"ahler }
\def\Kno{K\"ahler}

\def\KE{K\"ahler-Einstein }  \def\KEno{K\"ahler-Einstein}

\def\po1{partition of unity }

\def\proof{\hglue-\parindent{\it Proof. }}

\def\Ctwobefoot{\calfootC^{2,\beta}}

\def\Cinf{\h{\cal C}^\infty}

\def\Linf{L^\infty}

\def\CinfM{\Cinf(M)}

\def\n#1#2{||#1||_{#2}}

\def\strutdepth{\dp\strutbox}
\def\specialstar{\vtop to \strutdepth{
    \baselineskip\strutdepth
    \vss\llap{$\star$\ \ \ \ \ \ \ \ \  }\null}}
\def\marginalstar{\strut\vadjust{\kern-\strutdepth\specialstar}}
\def\marginal#1{\strut\vadjust{\kern-\strutdepth
    {\vtop to \strutdepth{
    \baselineskip\strutdepth
    \vss\llap{{ \small #1 }}\null} 
    }}
    }

\newcount\subsectionitemnumber
\def\clearsubsectionitemnumber{\subsectionitemnumber=0\relax}

\newcount\subsubsubsectionnumber
\def\clearsubsubsubsectionnumber{\subsubsubsectionnumber=0\relax}
\def\subsubsubsection#1{
\bigskip\noindent
\global\advance\subsubsubsectionnumber by 1%
{\rm
 \the\subsectionnumber.\the\subsubsectionnumber.\the\subsubsubsectionnumber}
{
#1.}
}

\newcount\subsubsectionnumber
\def\clearsubsubsectionnumber{\subsubsectionnumber=0\relax}
\def\subsubsection#1{
\clearsubsubsubsectionnumber
\bigskip\noindent
\global\advance\subsubsectionnumber by 1%
{%
\it \the\subsectionnumber.\the\subsubsectionnumber}
{
\it #1.}
}
\newcount\subsectionnumber
\def\clearsubsectionnumber{\subsectionnumber=0\relax}
\def\subsection#1{
\clearsubsectionitemnumber
\clearsubsubsectionnumber
\medskip\medskip\smallskip\noindent \global\advance\subsectionnumber by 1%
{%
\bf \the\subsectionnumber} 
{
\bf #1.}
}
\newcount\sectionnumber
\def\clearsectionnumber{\sectionnumber=0\relax}
\def\section#1{
\clearsubsectionnumber
\bigskip\bigskip\noindent \global\advance\sectionnumber by 1%
{%
\two
\the\sectionnumber} 
{
\two #1.}
}

\clearsectionnumber
\clearsubsectionnumber
\clearsubsubsectionnumber
\clearsubsubsubsectionnumber
\clearsubsectionitemnumber

\newcount\itemnumber
\def\clearitemnumber{\itemnumber=0\relax}
\def\c#1{ {\noindent\bf \the\itemnumber.} p. $#1$ \global\advance\itemnumber by 1}
\def\cn{ {\noindent\bf \the\itemnumber.} \global\advance\itemnumber by 1}

\clearitemnumber

\def\subsectionno#1{
\medskip\medskip\smallskip\noindent%
{%
\bf #1.}
}

\def\last#1{\advance\eqncount by -#1(\the\eqncount)\advance\eqncount by #1}
\def\llast{\advance\eqncount by -1(\the\eqncount)\advance\eqncount by 1}
\def\lllast{\advance\eqncount by -2(\the\eqncount)\advance\eqncount by 2}
\def\llllast{\advance\eqncount by -3(\the\eqncount)\advance\eqncount by 3}

\newcount\notenumber
\def\clearnotenumber{\notenumber=0\relax}
\def\note#1{\advance\notenumber by 1\footnote{${}^{(\the\notenumber)}$}
  {\lineskip0pt\notsosmall #1}}
\def\notewithcomma#1{\advance\notenumber by
1\footnote{${}^{\the\notenumber}$}
  {\lineskip0pt\notsosmall #1}}
                                                                                               
\clearnotenumber

\def\putnumber{%
\global\advance\subsectionitemnumber by 1{\the\subsectionnumber}.{\the\subsectionitemnumber}}
\def\numbering{{{\the\subsectionnumber}.{\the\subsectionitemnumber}}}

\def\FThm#1{\bigskip
{\noindent%
{\bf {T}heorem }
{\bf \putnumber.} {\it #1}\bigskip 
}}

\def\FCor#1{\bigskip
{\noindent%
{\bf {C}orollary }
{\bf \putnumber.} {\it #1}\bigskip
}}
\def\FLem#1{\bigskip
{\noindent%
{\bf {L}emma }
{\bf \putnumber.} {\it #1}\bigskip
}}

\def\FProp#1{\bigskip
{\noindent%
{\bf {P}roposition }
{\bf \putnumber.} {\it #1}\bigskip
}}

\def\FRemm#1{\medskip
{\noindent%
{\bf {R}emark }
{\bf \putnumber.}  {#1}\medskip
}}







\def\Abstract#1{
{\narrower\bigskip\bigskip
\noindent {{\bf Abstract.\ \ } #1}

}
}

\def\ref#1{{\bf[}{\sans #1}{\bf]}}

\def\reffoot#1{{\bfeight[}{\sansfoot #1}{\bfeight]}}

\def\opcit{\underbar{\phantom{aaaaa}}}

\def\sm{\smallskip}

\def\boxit#1{\vbox{\hrule\hbox{\vrule\kern3pt\vbox{\kern3pt#1\vglue3pt
\kern3pt}\kern3pt\vrule}\hrule}}

\long\def\frame#1#2#3#4{\hbox{\vbox{\hrule height#1pt
 \hbox{\vrule width#1pt\kern #2pt
 \vbox{\kern #2pt
 \vbox{\hsize #3\noindent #4}
\kern#2pt}
 \kern#2pt\vrule width #1pt}
 \hrule height0pt depth#1pt}}}

\def\kou{\frame{.3}{.5}{1pt}{\phantom{a}}}

\def\help{\ifmmode\aftergroup\noindent\quad\else\quad\fi}
\def\helpp{\ifmmode\aftergroup\noindent\hfill\else\hfill\fi}
\def\done{\helpp\kou\medskip}

\def\Ho{{\calH_{\o}}}

\def\on{\o^n}
\def\onminus1{\o^{n-1}/(n-1)!}
\def\ovp{\o_{\vp}}

\def\vpt{\vp_{\!t}}

\def\ovpta{\o_{\vp_\tau}}

\def\ovpt{\o_{\vpt}}

\def\ovpn{\o_{\vp}^n}
\def\ovptn{\o_{\vpt}^n}

\font\call=cmsy10   \font\Bbb=msbm10
\def\outlin#1{\hbox{\Bbb #1}}
\def\H{\hbox{\call H}}  \def\Ho{\hbox{\call H}_\omega}
 
\def\Cinf{{\call C}^\infty}  \def\RR{{\outlin R}}
\def\frac#1#2{{{#1}\over{#2}}}
\newcount \eqncount
\def \eqnno{\global \advance \eqncount by 1 \futurelet \nexttok \parsenexttok}
\def \eqn{\global \advance \eqncount by 1 \eqno\futurelet \nexttok \parsenexttok}
\def \eqnd{\global\advance \eqncount by 1 \futurelet\nexttok\parsenexttokd}
\def \parsenexttok{\ifx \nexttok $\Nomark\else\expandafter \Mark\fi}
\def \parsenexttokd{\ifx \nexttok \hfil\Nomark\else\expandafter \Mark\fi}
\def \Nomark {(\the \eqncount)}\def \Mark #1{\xdef #1{(\the \eqncount)}#1}

\let\exa\expandafter
\catcode`\@=11 
\def\CrossWord#1#2#3{%
\def\@x@{}\def\@y@{#2}
\ifx\@x@\@y@ \def\@z@{#3}\else \def\@z@{#2}\fi
\exa\edef\csname cw@#1\endcsname{\@z@}}
\openin15=\jobname.ref
\ifeof15 \immediate\write16{No file \jobname.ref}%
   \else \input \jobname.ref \fi 
\closein15
\newwrite\refout
\openout\refout=\jobname.ref 
\def\warning#1{\immediate\write16{1.\the\inputlineno -- warning --#1}}
\def\Ref#1{%
\exa \ifx \csname cw@#1\endcsname \relax
\warning{\string\Ref\string{\string#1\string}?}%
    \hbox{$???$}%
\else \csname cw@#1\endcsname \fi}
\def\Tag#1#2{\begingroup
\edef\head{\string\CrossWord{#1}{#2}}%
\def\writeref{\write\refout}%
\exa \exa \exa
\writeref \exa{\head{\the\pageno}}%
\endgroup}

\catcode`\@=12

\def\Tagg#1#2{\Tag{#1}{#2}}  

\def\TaggThm#1{\Tagg{#1}{Theorem~\numbering}}
\def\TaggLemm#1{\Tagg{#1}{Lemma~\numbering}}
\def\TaggRemm#1{\Tagg{#1}{Remark~\numbering}}
\def\TaggCor#1{\Tagg{#1}{Corollary~\numbering}}
\def\TaggProp#1{\Tagg{#1}{Proposition~\numbering}}
\def\TaggSection#1{\Tagg{#1}{Section~\the\subsectionnumber}}
\def\TaggS#1{\Tagg{#1}{\S~\the\subsectionnumber}}
\def\TaggSubS#1{\Tagg{#1}{\hbox{\S\S\unskip\the\subsectionnumber.\the\subsubsectionnumber}}}
\def\TaggSubsection#1{\Tagg{#1}{Subsection~\the\subsectionnumber.\the\subsubsectionnumber}}
\def\TaggEq#1{\Tagg{#1}{(\the\eqncount)}}
\def\Taggf#1{ \Tagg{#1}{{\bf[}{\sans #1}{\bf]}} }

\def\JJJ{\h{\rm J}}
\def\Ric{\hbox{\rm Ric}\,}
\def\Ricnotsosmall{\hbox{\notsosmall Ric}\,}

\def\Hc{\H_{c_1}}
\def\Dc{\calD_{c_1}}
\def\Hcfoot{\Hfoot_{c_1}}
\def\HcG{\H_{c_1}(G)}

\def\HcGplus{\H^{+}_{c_1}(G)}
\def\Hcplus{\H^{+}_{c_1}}

\def\Hcfootplus{\Hfoot^{+}_{c_1}}

\def\O{\Omega}

\def\Rico{\Ric\o}
\def\Ricovp{\Ric\ovp}
\def\Ricovpt{\Ric\ovpt}

\def\Ho{\H_{\Omega}}
\def\Do{\calD_{\Omega}}
\def\HoG{\H_{\Omega}(G)}
\def\HoGplus{\H^{+}_{\O}(G)}
\def\Hoplus{\H^{+}_{\O}}
\def\Mo{\calM_{\o}}

\def\MA{Monge-Amp\`ere }
\def\oone{\o_{\vp_1}}
\def\otwo{\o_{\vp_2}}

\def\Vm{V^{-1}}  \def\V{\frac1V}
\def\vpta{\vp_{\tau}}

\def\fovp{f_{\ovp}}

\def\Dvp{\Delta_\vp}

\magnification=1100
\hoffset1.1cm
\voffset1.3cm
\hsize5.09in
\vsize6.91in

\headline={\ifnum\pageno>1{\ifodd\pageno \oddheadline\else\evenheadline\fi}\fi}
\def\oddheadline{\centerline{\caps Energy functionals and \KE metrics}}
\def\evenheadline{\centerline{\caps Y. A. Rubinstein}}

\overfullrule0pt
\parindent12pt

\noindent
\centerline{\twosans On energy functionals, \KE metrics, and the}
\medskip
\centerline{\twosans Moser-Trudinger-Onofri neighborhood}
\smallskip

\font\tteight=cmtt8

\bigskip
\centerline{\sans Yanir A. Rubinstein%
\footnote{$^*$}{\cmreight Massachusetts Institute of Technology. 
Email: {\tteight yanir@member.ams.org}
\hfill\break
\vglue-0.5cm
\hglue-\parindent\cmreight
Current address: Department of Mathematics, Princeton University, Princeton, NJ 08544.}}

\def\centereps#1#2#3{\vglue#2\relax\centerline{\hbox to#1%
            {\special{eps:#3.eps x=#1 y=#2}\hfil}}}
\def\Hconeplus{\Hcplus(\Lambda_1)}
\font\cmreight=cmr8

\vglue0.03in

\footnote{}{\hglue-\parindent\cmreight December \thhnotsosmall{15}, 2006.
Revised October 2007.\hfill\break 
\vglue-0.5cm
\hglue-\parindent\cmreight
\vglue-0.5cm 
\hglue-\parindent\cmreight
32W20, %
53C25, %
58E11.} %
    
\vglue-0.12in

\bigskip
\Abstract{We prove that the existence of a \KE metric on a Fano manifold is 
equivalent to the properness of the energy functionals defined by Bando, Chen,
Ding, Mabuchi and Tian on the set of \K metrics with positive Ricci curvature. We also
prove that these energy functionals are bounded from below on this set if and only 
if one of them is. This answers two questions raised by X.-X. Chen. 
As an application, we obtain a new proof of the classical Moser-Trudinger-Onofri inequality on the two-sphere,
as well as describe a canonical enlargement of the space of \K potentials on which this inequality holds
on higher-dimensional Fano \KE manifolds.
}

\bigskip

 
\bigskip

\medskip
\subsectionnumber=0\relax

\subsection{Introduction}
\TaggSection{SectionIntroduction}
Our main purpose in this article is to give a new analytic
characterization of \KE manifolds in terms of certain functionals
defined on the infinite-dimensional space of \K forms.
As a corollary of our approach we also obtain a new proof of the classical Moser-Trudinger-Onofri
inequality on the two-sphere as well as an optimal extension of it to higher-dimensional
Fano \KE manifolds.

A necessary condition for a manifold to admit a \KE metric 
is that its
first Chern class be either positive, negative or zero. 
Aubin and Yau proved that this condition is also sufficient in the second 
case and Yau proved that the same is true also in the third case.

Yet additional geometric assumptions are necessary in the first case (in this
case the manifold is called Fano): Matsushima proved that the group of automorphisms
must be reductive, Futaki proved that a certain character on the algebra of holomorphic
vector fields must be trivial, and Kobayashi and L\"ubke proved that the
tangent bundle must be stable.
Since then
much work has been done on the subject
(see for example the recent expositions \ref{Bi,F2,Sz,Th}).

In this article we will restrict attention to two closely related analytic criteria 
relating the existence of \KE metrics to properties of certain 
energy functionals (see the end of this section and \Ref{SectionEnergyFunctionals} for notation and definitions) 
on the space of \K forms $\Hc$. The first, introduced by Tian, can be thought of as a 
``stability" criterion \ref{T3}. It expresses
the existence of a \KE metric as equivalent to the properness of an energy functional:

\FThm{\TaggThm{TianThm}\Tagg{TianThmno}{\numbering}
\ref{T3,T4,TZ}  Let $(M,\JJJ)$ be a Fano manifold and assume that 
$\Aut(M,\JJJ)$ is finite. 
Then the following are equivalent:
(i) $(M,\JJJ)$ admits a \KE metric,
(ii) $E_0$ is proper on $\Hcplus$,
(iii) $F$ is proper on $\Hcplus$.
}

\noindent
The finiteness assumption\footnote{$^1$}{\cmreight 
Since automorphisms of the complex structure preserve the first Chern class this assumption
is\break
\vglue-0.5cm
\hglue-\parindent\cmreight
 equivalent to the triviality of $\mf\hbox{\cmreight aut}(M,\hbox{\cmreight J})$ \reffoot{Fu, Theorem 4.8}.} 
 covers, for example, all \KE Fano surfaces 
except the product of two Riemann spheres, the projective plane $\PP^2$, and 
$\PP^2$ blown up at 3 non-collinear points \ref{S1,T2,TY}.
However, there is a slightly more technically involved version of \Ref{TianThm}
\unskip, also due to Tian, which applies to
all \KE Fano manifolds, that will be stated in \Ref{SectionProofs} 
(\Ref{TianFullThm}).

The second analytic criterion, introduced by Bando and Mabuchi, can be thought 
of as a ``semi-stability" condition \ref{BM}. Two related formulations appeared subsequently \ref{Ba,DT}. 
It expresses the existence of ``almost" \KE metrics
as a consequence
of the lower boundedness of an energy functional:

\FThm{\TaggThm{BDTThm}\Tagg{BDTThmno}{\numbering}
\ref{BM,DT}
Let $(M,\JJJ)$ be a Fano manifold.
Assume that either $F$ or $E_0$ is bounded from below on $\Hcplus$ and let $\eps>0$. 
Then $(M,\JJJ)$ admits a 
\K metric $\o_\eps\in\Hc$ satisfying
$\Ric\o_\eps>(1-\eps)\o_\eps$.
}

It is worth mentioning that a precise characterization of Fano manifolds for which these functionals
are bounded from below is still lacking. Also, examples of such manifolds which
are not \KE are yet to be given.

We point out that \Ref{TianThm} and the version of \Ref{BDTThm} for the functional $F$ were originally stated
with the assumptions on properness and boundedness made on the whole space 
of \K forms $\Hc$ rather than on the subspace
of forms of positive Ricci curvature $\Hcplus$. 
However, the respective existence proofs only make use of those assumptions on $\Hcplus$. 
Thus \Ref{TianThm} implies that the properness of the functionals on $\Hcplus$ implies 
their properness on $\Hc$.
In addition, in \Ref{RemmBound} we prove that for any Fano manifold also the lower boundedness of the 
functionals on $\Hcplus$ implies their lower boundedness on $\Hc$. Therefore it seems more natural to state
Theorems \Ref{TianThmno} and \Ref{BDTThmno} in the equivalent manner above. 
This will also be justified by the results of \Ref{SectionMTO} (in particular 
\Ref{EnunboundedCor}).

Chen and Tian constructed a family of energy functionals $E_1,\ldots,E_n$, analogues of the 
`K-energy' (`\K energy') $E_0$ corresponding to higher degree elementary symmetric polynomial expressions of the eigenvalues
of the Ricci tensor \ref{CT}. As with $E_0$ and $F$, \KE metrics are critical points of these 
functionals  
and it is therefore a natural idea to seek to extend Theorems \Ref{TianThmno} 
and \Ref{BDTThmno} to $k=1,\ldots,n$.
In this direction, an analogue of \Ref{TianThm} for $k=1$ was proved recently 
by Song and Weinkove \ref{SW}. 
The main purpose of the present article is to prove the following two statements:

\FThm{\TaggThm{PropernessThm}\Tagg{PropernessThmno}{\numbering}
Let $(M,\JJJ)$ be a Fano manifold and assume that $\Aut(M,\JJJ)$ is finite. 
Let $k\in\{0,\ldots,n\}$. Then the following are equivalent:
(i) $(M,\JJJ)$ admits a \KE metric,
(ii) $E_k$ is proper on $\Hcplus$,
(iii) $F$ is proper on $\Hcplus$.
}

\vskip-0.3cm

\FThm{\TaggThm{BoundednessThm}\Tagg{BoundednessThmno}{\numbering}
Let $(M,\JJJ)$ be a Fano manifold
and let $k\in\{0,\ldots,n\}$. Assume that either $F$ or $E_k$ is bounded from below on $\Hcplus$ and let $\eps>0$. 
Then $(M,\JJJ)$ admits a  
\K metric $\o_\eps\in\Hc$ satisfying $\Ric\o_\eps>(1-\eps)\o_\eps$.
}

Our proofs carry over to \KE manifolds admitting holomorphic vector fields (for the 
more general statements the reader is referred to 
Sections \Ref{SectionContinuityno} and \Ref{SectionProofsno}). 
We remark that while 
\Ref{PropernessThm} generalizes the work of Song and Weinkove, our
methods provide a considerable simplification over the ones used there.

These theorems show that the functionals $E_k$ are, on the one hand,
closely related to geometric stability, and, on the other hand, all equivalent
in a suitable sense.\footnote{$^2$}{\cmreight In particular, 
after posting the first version of this article I became 
aware of the fact that Theo-\break
\vglue-0.5cm
\hglue-\parindent\cmreight
rems
\Ref{PropernessThmno} and
\Ref{BoundednessThmno}
 answer questions posed recently by Chen \reffoot{CLW}.}

To prove these theorems we first observe that a certain formula of Bando and Mabuchi
for the `Ricci energy' $E_n$ extends naturally to all of the functionals $E_k$. 
The merit of this new formula (\Ref{BMRProp}) is that it succinctly captures the relation between
the different functionals.
This shows in particular that
the lower boundedness of $E_k$ implies the lower boundedness of $E_{k+1}$. We then interpret another observation 
of Bando and Mabuchi in order to close the loop and prove that the lower boundedness of $E_n$ on $\Hcplus$
implies that of $F$ on $\Hc$. This step is crucial in proving \Ref{BoundednessThm}. In fact it proves more, namely,
that the lower boundedness of any one of the functionals implies that of the rest 
(\Ref{CorrBound}). Special cases of this fact have been observed previously \ref{CLW,DT,L2,P}
(see \Ref{RemmBound}).

To prove \Ref{PropernessThm} we consider the continuity method path
\Ref{MASecondpathseq}
introduced by Aubin \Ref{A2}.
As before, we show that the properness of $E_k$ implies the properness of $E_{k+1}$.
Next, assuming $E_n$ is proper and using \Ref{BoundednessThm} we conclude that this 
path exists for all $t\in[0,1)$.
We show that on a fixed 
interval $[t_0,1)$ each of the functionals $E_k$ is uniformly bounded from above with $t_0$
depending only on $n$, and then conclude.

In \Ref{SectionMTO} we observe that \Ref{BMRProp} allows to obtain
without additional effort 
a strengthened version of 
the second main result of Song and Weinkove, the one concerning the nonnegativeness 
of the energy functionals with respect to a \KE base metric. 
We observe that our results, when combined with previous ones
\ref{BM,DT}, provide for a new and 
entirely \Kahler geometric
proof of the Moser-Trudinger-Onofri inequality on the Riemann sphere. As a corollary of 
this approach we also characterize the functions for which this inequality continues
to hold on higher-dimensional \KE Fano manifolds, thus extending the work of Ding and Tian. 
We call the set of all such functions the Moser-Trudinger-Onofri neighborhood of the space of
\K potentials. It is a canonically defined set that strictly contains the space of \K potentials and lies
within $\CinfM$. 
This provides a higher dimensional analogue of the original Moser-Trudinger-Onofri inequality
that is optimal in a certain sense and brings Ricci curvature into the picture (\Ref{CharacterizationThm}).
Finally, we are able to show that the energy functionals $E_2,E_3,\ldots$ are not bounded
from below on $\Hc$ (\Ref{EnunboundedCor}).

The results herein have applications also to the study of the \Kno-Ricci flow and 
geometric stability \ref{R} that will appear in a subsequent article.

The article is organized as follows. In \Ref{SectionEnergyFunctionals} we review the relevant background concerning energy functionals
and present the formula for the functionals $E_k$ (\Ref{BMRProp}) whose proof appears
in the Appendix. 
In \Ref{SectionContinuity} we review results concerning the continuity method approach.
The proofs of our main results
are contained in \Ref{SectionProofs}. \Ref{SectionMTO} concludes with our results 
on the lower boundedness of the functionals $E_k$ and on the generalized Moser-Trudinger-Onofri inequality.

\subsectionno{Setup and notation}
Let $(M,\JJJ)$ be a connected compact closed \K manifold of complex dimension $n$
and let $\O\in H^2(M,\RR)\cap H^{1,1}(M,\CC)$ be a \K class
with $d=\del+\dbar$. Define the Laplacian $\D=-\dbar\circ\dbar^\star-\dbar^\star\circ\dbar$ with
respect to 
a Riemannian metric $g$ on $M$ and assume that $\JJJ$ is compatible with $g$ and parallel
with respect to its Levi-Civita connection. 
Let
$\gHerm=1/\pi\cdot\gij(z)\dz i\otimes\dzb j$ be the associated \K metric, 
that is the induced Hermitian metric on $(T^{1,0}M,\JJJ)$, and let
$\o:=\o_g=\i/2\pi\cdot\gij(z)\dzidzjb$ denote its corresponding \K form,
a closed positive $(1,1)$-form on $(M,\h{\rm J})$
such that $\gHerm=\ha g-\frac\i2\o$. Similarly denote by $g_\o$ the Riemannian
metric induced from $\o$ by $g_\o(\cdot,\cdot)=\o(\cdot,\JJJ\,\cdot)$. 
For any \K
form we let $\Ric(\omega)=-\i/2\pi\cdot\ddbar\log\det(\gij)$ denote the
Ricci form of $\o$. It is well-defined globally and represents
the first Chern class $c_1:=c_1(T^{1,0}M,\h{\rm J})\in H^2(M,\ZZ)\cap H^{1,1}(M,\CC)$.
One calls $\o$ \KE if $\Ric\o=a\o$ for some real $a$.

Denote by $\Do$
the space of 
all closed $(1,1)$-forms whose cohomology class is $\O$.
For a \K form $\o$ with $[\o]=\O$ we will consider the space of strictly $\o$-plurisubharmonic functions
$$
\calH_\o=\{\vp\in\CinfM\,:\, \ovp:=\o+\i\ddbar\vp >0\},
$$
and the subspace $\Ho\sseq\Do$ of \K forms cohomologous to $\O$. We denote by $\Ho^+\sseq\Ho$ the 
subspace of those \K forms whose Ricci curvature is positive.
Let $\Aut(M,\JJJ)$ denote the complex Lie group of automorphisms (biholomorphisms)
of $(M,\JJJ)$
and denote by  $\aut(M,\JJJ)$ its Lie algebra of infinitesimal automorphisms composed
of real vector fields $X$ satisfying $\L_X\JJJ=0$. 
Let $G$ be 
any compact real Lie subgroup of $\Aut(M,\JJJ)$, and 
let $\Aut(M,\JJJ)_0$ denote the identity component
of $\Aut(M,\JJJ)$.
We denote by $\HoG\sseq\Ho$ and $\HoGplus\sseq\Hoplus$ the corresponding subspaces of $G$-invariant forms.

\subsection{Certain energy functionals on the space of \K forms}
\TaggSection{SectionEnergyFunctionals}
We call a real-valued function $A$ defined on a subset $\Dom(A)$ of 
$\Do\times\Do$ an energy functional if it 
is zero on the diagonal restricted to $\Dom(A)$. By
a Donaldson-type functional, or exact energy functional, we will mean an energy 
functional
that satisfies the cocycle condition $A(\o_1,\o_2)+A(\o_2,\o_3)=A(\o_1,\o_3)$ with each
of the pairs appearing in the formula belonging to $\Dom(A)$
\ref{Do,M,T4}. We will occasionally refer to both of these simply as 
functionals and exact functionals, respectively.
Note that if an exact functional is defined on $U\times W$ with $U\sseq W$ then there exists a unique exact
functional defined on $W\times W$ extending it.

Let $V:=\int_M\o^n=[\o]^n([M])$.
The energy functionals
$I,J$, introduced by Aubin \ref{A2}, are defined for each pair 
$(\o,\o_{\vp}:=\o+\i\ddbar\vp)\in\Do\times\Do$ by
$$\eqano{
I(\o,\ovp) & =\Vm\int_M\i\del\vp\w\dbar\vp\w\sum_{l=0}^{n-1}\o^{n-1-l}\w\ovp^{l}
=\Vm\intM\vp(\on-\ovpn),&\eqnno\Ieq\cr
J(\o,\ovp) & =\frac{\Vm}{n+1}\int_M\i\del\vp\w\dbar\vp\w\sum_{l=0}^{n-1}(n-l)\o^{n-l-1}\w\ovp^{l}.&\eqnno\Jeq\cr
}$$
One may also define
them via a variational formula. 
Connect each pair
$(\o,\o_{\vp_1}:=\o+\i\ddbar\vp_1)$ with a piecewise smooth path
$\{\ovpt\}_{t\in[0,1]}$ (we regard this path as a function
on $M\times[0,1]$ and occasionally suppress the subscript $t$). Then we have for any such path
$$
\eqano{
(I-J)(\o,\o_{\vp_1}) & =
-\V\int_{M\times [0,1]}\vpt n\i\ddbar\dot\vpt\w\ovpt^{n-1}\w dt,&\eqnno\IminusJeq
\cr
J(\o,\o_{\vp_1}) & =\V\int_{M\times [0,1]}\dot\vpt(\on-\ovptn)\w dt.&\eqnno\Jeq
}$$
On $\Ho\times\Ho$ $I,J$ and $I-J$ are all nonnegative (and hence non-exact) and equivalent, namely,
$$ 
\frac1{n^2}(I-J)\le \frac1{n(n+1)}I\le\frac1n J\le I-J\le \frac{n}{n+1} I\le nJ.
\eqn\IJineq
$$
Note that pulling-back both arguments of these functionals by an automorphism of $(M,\h{\rm J})$ does not change their value. It is important to understand the behavior of these functionals
also outside the subspace $\Ho$:
\FLem{
\TaggLemm{UnboundedLemma}
Let $\o\in\Ho$. Then $I(\o,\cdot)$ is unbounded from above on $\Ho$ and, when $n>1$, unbounded
on $\Do$.
}
\proof Fix a holomorphic coordinate patch   
$$
\psi:U\ra\CC^n,\;
\psi(q)={\bf z}(q):=(z^1(q),\ldots,z^n(q)), \;\all\,q\in U\sseq M.
$$
Let $a>0$ be such that $\ii{\psi}(\{v\in\CC^n\,:\, |v|<3a\})\sseq U$.
For the first statement, define $\tilde\vp_{b}$ 
by letting $\tilde\vp_{b}=b|{\bf z}|^2$ on 
$\ii{\psi}(\{v\in\CC^n\,:\, a<|v|<2a\})$ and constant elsewhere on $U$
in such a way that it is continuous. Approximate $\tilde\vp_{b}$ by
smooth functions $\vp_{b,m}$ that agree with it outside the set
$\ii\psi\{\{v\in\CC^n\,:\, |v|\in(a-\frac1m,a+\frac1m)\cup(2a-\frac1m,2a+\frac1m)\})$
and that satisfy $|\vp_{b}-\vp_{b,m}|<\frac1m$ on $U$. Given $a_2>0$
there exists $b$ and a corresponding $m$ such that 
$\vp_{b,m}\in\Mo$ and 
$I(\o,\o_{\vp_{b,m}})>a_2$.

For the second statement, construct similarly functions, as above,
now setting $\tilde\vp_{b}=-b(|z_1|^2+|z_2|^2)$
on $\ii{\psi}(\{v\in\CC^n\,:\, a<|v|<2a\})$. Again one may approximate using
functions $\vp_{b,m}$. Expanding $(\o+\i\ddbar\vp)^l$ using the binomial formula
it then follows that up to a term that is uniformly bounded for $m$ sufficiently large,
$I(\o,\o_{\vp_{b,m}})$ equals
$\Vm\int_M\i\del\vp\w\dbar\vp\w\o^{n-2}\w(a_2\o+a_3\i\ddbar\vp_{b,m})$
for some $a_2,a_3>0$.
We then see that given any $a_4>0$ there exists $b$ and a corresponding $m$ such
that $I(\o,\o_{\vp_{b,m}})<-a_4$.\done

We say that an exact functional $A$ is bounded from below 
on $U\sseq\Ho$ if for every $\o$ such that $(\o,\ovp)\in\Dom(A)$ and $\ovp\in U$ holds $A(\o,\ovp)\ge C_\o$
with $C_\o$ independent of $\ovp$. 
We say it is proper (in the sense of Tian) on a set $U\sseq\HoG$  if for each $\o\in\HoG$ 
there exists a smooth function 
$\tau_\o:\RR\ra\RR$ 
satisfying $\lim_{s\ra\infty}\tau_\o(s)=\infty$ such that 
$A(\o,\ovp)\ge\tau_\o((I-J)(\o,\ovp))$
for every $\ovp\in U$. This is well-defined, in other words
depends only on $[\o]$ since the failure of $I-J$ to satisfy the cocycle condition is 
under control with respect to the two base metrics, $\o$, $\oone$ say, to wit,
$$
(I-J)(\o,\otwo)-(I-J)(\oone,\otwo)=
(I-J)(\o,\oone)-\V\intM\vp_1(\otwo^n-\oone^n)
,
$$  
with the last term controlled by the oscillation of $\vp_1$. 
Properness of a functional implies it has a lower bound.

Define the following collection of energy functionals for each $k\in\{0,\ldots,n\}$
$$\eqano{
I_k(\o,\ovp)
= &\
\V\intm\i\del \vp\w\dbar \vp\w\sum_{l=0}^{k-1} \frac{k-l}{k+1}\o^{n-1-l}\w\ovp^l
\cr
= &\ 
\frac\Vm{k+1}\intm\vp(k\on-\sum_{l=1}^k\o^{n-l}\w\ovp^l). &\eqnno\Ikeq
}$$
Note that $I_n=J,\; I_{n-1}=({(n+1)J-I})/n$. 

Chen and Tian \ref{CT} defined another such family 
$$
J_k(\o,\o_{\vp_1})
=
\Vm\int_{M\times [0,1]}
\dot\vpt(\ovpt^{k}\w\o^{n-k}-\ovptn)\w dt,\q k=0,\ldots,n. \eqn\JkVariationeq
$$
(Note that $J_{n-k-1}/(k+1)$ in their article corresponds to $J_k$ in this article.)
The following computation relates these two families of functionals.

\FLem{
\TaggLemm{IkJkLemma} The following relation holds on $\Ho\times\Ho$:
$$
I_k(\o,\ovp)=
J(\o,\ovp)-J_{k}(\o,\ovp).
$$
}
\proof Given a path $\{\ovpt\}_{t\in[0,1]}$ 
we compute the variational equation for $I_k$.
$$
\eqano{
(k+1)\ddt I_k(\o,\ovpt)
= & 
-\V\intm\sum_{l=0}^{k-1}\Big(
2\dot\vp\i\ddbar\vp\w (k-l)\o^{n-1-l}\w\ovp^l
\cr
& +\vp\i\ddbar\vp\w\i\ddbar\dot\vp\w l(k-l) \o^{n-1-l}\w\ovp^{l-1}\Big)
\cr
= & 
-\V\intm\dot\vp\i\ddbar\vp\w\sum_{l=0}^{k-1}\Big(
2 (k-l)\o^{n-1-l}\w\ovp^l
\cr
& +(\ovp-\o)\w l(k-l) \o^{n-1-l}\w\ovp^{l-1}\Big)
\cr
= & 
-\V\intm\dot\vp\i\ddbar\vp\w\Big(
\sum_{l=0}^{k-1}2 (k-l)\o^{n-1-l}\w\ovp^l
\cr
& +\sum_{l=1}^{k-1} l(k-l) \o^{n-1-l}\w\ovp^{l}
\cr
&-\sum_{l=0}^{k-2} (k-l-1)(l+1)\o^{n-1-l}\w\ovp^{l}
\Big)
\cr
= & 
-(k+1)\V\intm\dot\vp\i\ddbar\vp\w
\sum_{l=0}^{k-1}\o^{n-1-l}\w\ovp^l,
}$$
and putting $\i\ddbar\vp=\ovp-\o$ we have
$$
\ddt I_k(\o,\ovpt)=
\Vm\intm\dot\vpt
(\on-\o^{n-k}\w\ovpt^k).\eqn\IkVariationeq
$$
Combining with \JkVariationeq\ and \Jeq\ we conclude.\done
Note that from the definitions it follows that
$$
0\le I_{k}(\o,\ovp)\le J(\o,\ovp),\qq \h{on\ \ } \Ho\times\Ho.\eqn\IkJeq
$$
As a corollary
of \Ref{IkJkLemma} we have therefore
$0\le J_{k}(\o,\ovp)\le J(\o,\ovp)$ on $\Ho\times\Ho$. We point out
that this upper bound improves \ref{CT,Corollary 4.5} while the lower
bound appears to be new. Also from \Ikeq
$$
\frac{k+2}{k+1}I_{k+1}\ge \frac{k+1}{k}I_{k}, 
\qq \h{on\ \ } \Ho\times\Ho. \eqn\Ikineq
$$ 
Note that in particular $I_{k+1}\ge I_k$ and
so by \Ref{IkJkLemma}
$J_{k}\ge J_{k+1}$.
We note in passing that this lemma also yields the following formula
$$
\eqano{
J_k(\o,\ovp)=\frac\Vm{n+1}\intm\i\del\vp\w\dbar\vp\w
\Big( &
\frac{n-k}{k+1}\sum_{l=0}^{k-1}(l+1)\o^{n-1-l}\w\ovp^l
\cr
& +\sum_{l=k}^{n-1}(n-l)\o^{n-1-l}\w\ovp^l
\Big).
&\eqnno\Jkeq
}$$

\def\Riceight{\hbox{\cmreight Ric}\,}

The energy functionals $E_k,\,\,k=0,\ldots,n$, are
defined by 
$$\eqalignno{
E_k(\omega,\omega_{\vp_1}) \,\,=\,\,\,& \Vm
\int_{M\times [0,1]}
\D_{\vp_t}\dot\vp_t\Ric(\o_{\vp_t})^k\w\ovpt^{n-k}\w dt &\eqnno\Ekdefeq \cr 
&
-\frac{n-k}{k+1}\Vm \int_{M\times [0,1]}
\dot\vp_t(\Ric(\o_{\vp_t})^{k+1}-\mu_k\ovpt^{k+1})\w\ovpt^{n-1-k}\w
dt,
}$$ where
$\mu_k:=\frac{c_1^{k+1}\cup[\o]^{n-k-1}([M])}{[\o]^n([M])}.\,$ 
This gives rise to
well-defined exact energy functionals
\ref{CT} (note that $E_k/(k+1)$ in the aforementioned article corresponds to $E_k$ in this article).
The K-energy, $E_0$, was introduced by Mabuchi \ref{M}, while
$E_n$, which we refer to as the `Ricci energy', was introduced by Bando and Mabuchi%
\footnote{$^3$}{\cmreight 
\KE forms are the only critical points of these two functionals when 
$\mf \O=\mu c_1, \mu\in\{\pm1\}$:
\break
\vglue-0.5cm
\hglue-\parindent\cmreight
For $\mf E_0$ see \reffoot{T4, p. 19}
while for $\mf E_n$ the
critical forms satisfy $\mf (\mu\Riceight\o)^n=\o^n$ and 
writing $\mf \mu\Riceight\o=$
\break
\vglue-0.5cm
\hglue-\parindent\cmreight
$\mf\o+\i\ddbar f$ we see that $\mf \mu\Riceight\o>0$ at the
minimum of $\mf f$. Since
the smallest eigenvalue of a H\"older
\break
\vglue-0.5cm
\hglue-\parindent\cmreight
continuous matrix-valued function is also
H\"older
continuous \reffoot{Al, p. 438} we conclude that $\mf \mu\Riceight\o>0$ 
\break
\vglue-0.5cm
\hglue-\parindent\cmreight
implying
that $\mf f$ is constant by the uniqueness argument of Calabi
(for a 
different proof see \reffoot{Ma, \S8}). 
\break
\vglue-0.5cm
\hglue-\parindent\cmreight
However when $\mf c_1=0$ there are nontrivial solutions of 
$\mf (\Riceight\o)^n=0$ if the
manifold is a product.
For 
\break
\vglue-0.5cm
\hglue-\parindent\cmreight
$\mf \mu=1$ critical points of $\mf E_k$ with nonnegative 
Ricci curvature are necessarily \KE 
\reffoot{To}.
}
 \ref{BM}.

For each $\o\in\Ho$ these functionals (being exact) induce a (real) Lie group homomorphism
$\Aut(M,\JJJ)_0\ra\RR$ given by $h\mapsto E_k(\o,h^\star\o)$. The corresponding Lie
algebra homomorphism $\aut(M,\JJJ)\ra\RR$ is given by 
$X\mapsto\ddt\big|_0 E_k(\o,(\exp tX)^\star\o)$. This naturally extends to a 
complex Lie algebra homomorphism
$$X\mapsto
\calF_k(X;\o):=\ddt\Big|_0 E_k(\o,(\exp tX)^\star\o)
-\i
\ddt\Big|_0 E_k(\o,(\exp tJX)^\star\o).\eqn\Fkdefeq
$$
Changing $\o$ within a fixed cohomology class does not
change the homomorphism \ref{CT,Ma}. 
This is an extention of
the Bando-Calabi-Futaki Theorem, the case $k=0$ 
\ref{Be,C,F1}
(the construction was further generalized by Futaki \ref{F3}).
One calls these homomorphisms
Futaki characters (or invariants). 
When
$(M,\JJJ,\o)$ is Fano \KE it follows from \Ekdefeq\ that $\calF_k$ is trivial and hence
$E_k(\o_,\ovp)=0$ if $\ovp$ is \KEno, since the set
of \KE metrics is equal to an $\Aut(M,\JJJ)_0$-orbit of $\o$ \ref{BM}.

Unless otherwise stated, from now and on we will assume that $(M,\JJJ)$ is Fano 
and let $\ovp\in\Hc$. Let $f_{\ovp}\in\CinfM$ denote 
the unique function satisfying $\i\ddbar f_{\ovp}=\Ricovp-\ovp$ 
and $\Vm\intm e^{f_{\ovp}}\ovpn=1$.
Following Ding \ref{D}, define an exact functional on $\Hc\times\Dc$ by
$$
F(\omega,\omega_\vp)
 =
J(\o,\ovp)-\V\intm\vp\on-\log\V\intM e^{f_\o-\vp}\o^n.
$$
The critical points of this functional are the \KE metrics.
We state the following relation between the functionals $E_0$ and $F$.
\FLem{\TaggLemm{DTLemma} \ref{DT} Let $(\o,\ovp)\in\Hc\times\Hc$. Then
$$
F(\o,\ovp)=
E_0(\o,\ovp)
+\V\intM f_{\ovp}\ovpn
-\V\intM f_\o \on.
$$
}
Note that 
$$
\V\intm f_{\ovp}\ovpn\le 
\V\intm e^{f_{\ovp}}\ovpn-1=0.\eqn\fovpNegativeeq
$$

Note also that one may define a Lie algebra homomorphism corresponding to $F$ similarly to the
construction for $E_k$ in \Fkdefeq. \Ref{DTLemma} implies that this homomorphism will coincide
with $\calF_0$.

An equivalent form of the following was stated 
by Bando and Mabuchi \ref{BM,(1.5)}. 
\FLem{\TaggLemm{FzeroEn}
For every $(\o,\ovp)\in\Hcplus\times\Hc$ one has 
$$
E_n(\o,\ovp)=F(\Ric\o,\Ric\ovp).
$$
}
\noindent Note that by exactness this formula completely determines
$E_n$ on $\Hc\times\Hc$, as remarked earlier. 

\sm
\proof
Let $\{\vpt\}$ denote a smooth family of functions such that
$\o_{\vp_0}=\o,\; \o_{\vp_1}=\ovp$. 
Write $\Ric\ovpt=\Ric\o+\i\ddbar\log{\o^n\ov\ovptn}$. Then
$f_{\Ricnotsosmall\o}=\log{\on\ov(\Ricnotsosmall\o)^n}$. Thus for each $t\in[0,1]$,
$$
F(\Rico,\Ricovpt)
=
J(\Rico,\Ricovpt)-\V\intm\log{\o^n\ov\ovptn}(\Rico)^n.
$$
Hence,
$$
\eqa{
\ddt F(\Rico,\Ricovpt)
&=-\Vm\intM(-\D_t\dot\vpt)(\Ricovpt)^n
=\ddt E_n(\o,\ovpt),
}$$
from which we conclude by integration.\done

Bando and Mabuchi derived the following elegant formula.%
\FProp{\TaggProp{BMProp} 
\ref{BM, (1.8.1)}
For every $(\o,\ovp)\in\Hc\times\Hc$,
$$
E_n(\o,\ovp)=
E_0(\o,\ovp)+J(\ovp,
\Ricovp)-J(\o,
\Rico).
$$
}

We now show that \Ref{BMProp} can be generalized as follows.

\FProp{\TaggProp{BMRProp}
Let $k\in\{0,\ldots,n\}$. For every $(\o,\ovp)\in
\Hc\times\Hc$, 
$$\eqano{
E_k(\o,\ovp)
& =
E_n(\o,\ovp)-J_{k}(\ovp,
\Ricovp)+J_{k}(\o,
\Rico),&\eqnno\FirstBMRPropEq
\cr\cr & =
E_0(\o,\ovp)+I_k(\ovp,
\Ricovp)-I_k(\o,
\Rico),&\eqnno\SecondBMRPropEq
\cr\cr & =
((1-\hbox{$\frac{l}{k+1}$})E_0+\hbox{$\frac{l}{k+1}$}E_n)(\o,\ovp)+(I_k-\hbox{$\frac{l}{k+1}$}J)(\ovp,
\Ricovp)\qq&\eqnno\ThirdBMRPropEq
\cr\cr
& \qq-(I_k-\hbox{$\frac{l}{k+1}$}J)(\o,
\Rico),\q\all\,l\in\{0,\ldots,k+1\}.
}$$
}
The proof appears in the Appendix. 
The functionals $E_k$ are thus seen to be described as
`\Kno-Ricci' energies, ``interpolating" between the \K energy $E_0$
and the Ricci energy $E_n$.
We note that there exist counterparts of the formulas
presented in this section for other classes \ref{R}.

One particularly visible consequence of \Ref{BMRProp} is the fact that
the homomorphisms $\calF_k$ all coincide,
a result first proved by Maschler \ref{Ma, (17)} using an equivariant formulation and 
later by Liu \ref{Li, \S3} by a direct computation
(see also \ref{L1}). 
For other explicit expressions for the functionals $E_k$ see \ref{CT,L1,P,SW}.

\subsection{Continuity method approach}
\TaggSection{SectionContinuity}
\Tagg{SectionContinuityno}{\the\subsectionnumber}
Consider the path $\{\ovpt\}\sseq\Hc$ given implicitly by
$$\eqano{
\ovpt^n =  e^{(t+1)f_\o+c_t}\o^n,\qq & \q t\in[-1,0],\cr
\ovpt^n =  e^{f_\o-t\vpt}\o^n,\qq &  \q t\in[0,1],&\eqnno\MASecondpathseq\cr
\TaggEq{MASecondpathseq}
}$$
with the normalizations
$\int_Me^{(t+1)f_\o+c_t}\o^n=V$ for $t\in[-1,0]$ 
and  $\int_Me^{f_\o-t\vpt}\o^n=V$ for $t\in[0,1]$.
Note that the first segment always exists by the proof of the Calabi-Yau Theorem \ref{Y} while the 
second, when it exists, deforms the metric to a \KE metric \ref{A2}:
$$
\Ricovpt-\ovpt=-(1-t)\i\ddbar\vpt,\qq\q t\in[0,1].\eqn\AubPatheq
$$
We will make use of the following Proposition:
\FProp{
\TaggProp{AubDef}
\ref{BM, Theorem 5.7} Assume that $(M,\JJJ)$ is Fano and 
let $G$ be a compact subgroup of $\Aut(M,\JJJ)$. Assume that 
$E_0$ is bounded from below on $\HcGplus$ and let $\o\in\HcG$. 
Then 
\MASecondpathseq\ 
has a unique smooth solution for each $t\in[0,1)$.
}

\noindent
Note that by \Ref{DTLemma} and \fovpNegativeeq\ the same conclusion holds with $E_0$ 
replaced by $F$.  
In particular, \Ref{BDTThm} is a direct corollary of \Ref{AubDef} combined with
this observation (one obtains a version of \Ref{BDTThm} with 
the free choice of a subgroup $G$, although this, as opposed to the refinement of \Ref{TianThm} that will be
given in the next section, should not be considered as a gain in generality).
We also note that one of the important ingredients in the proof of \Ref{AubDef} is the 
fact that $(I-J)(\o,\cdot)$ is nondecreasing along the continuity path
\MASecondpathseq\
\ref{BM, Theorem 5.1; T1, p. 232; T4, Lemma 6.25}.

It is worth noting that Bando has shown that if $\o\in\HcG$ satisfies 
$\Ric\o>(1-\eps)\o,\; \eps>0,$ then ``flowing" it along the Ricci flow 
will produce another metric in $\HcG$ whose scalar curvature differs from $n$ by
at most a fixed constant times $\eps$ \ref{Ba}. Therefore, 
the existence of a lower
bound for $E_0$ or for $F$ implies the existence of \K metrics in $\HcG$ 
whose scalar curvature is as close to a constant as desired (the original
result of Bando
extends to the $G$-invariant setting since its proof
makes use of a \Kno-Ricci flow which, like the continuity method, preserves $\HcG$).
These can be thought of as ``almost \KEno" metrics since a \K metric of constant scalar
curvature in $\Hc$ is necessary \KEno.

\subsection{Boundedness and properness properties of energy functionals}
\TaggSection{SectionProofs}
\Tagg{SectionProofsno}{\the\subsectionnumber}
By Matsushima's 
Theorem, when a \KE form $\o$ exists 
the Lie algebra of Killing vector fields  
is a real form of $\aut(M,\JJJ)$ \ref{Be,Mat,S2}. %
In other words, when a \KE metric exists we may take $G$ to be the isometry group $\Iso(M,g_\o)$.
Also, $\aut(M,\JJJ)$ is then isomorphic to an eigenspace of the Laplacian, namely,
$$
\aut(M,\JJJ)
\isom
\Lambda_1:=\{\psi\in\CinfM\,:\, -\D_\o\psi=\psi\}.
$$
Set
$$
\Hconeplus:=\{\ovp\in\Hcplus\,:\, \intM\vp\psi\o^n=0,\q \all\psi\in\Lambda_1\}.
$$
Similarly, define $\Hc(\Lambda_1)$.
We may now state the following theorem of Tian which is a refined version of
\Ref{TianThm}.\footnote{$^3$}{\cmreight A detailed exposition of this theorem
will be found in a forthcoming article of Tian and Zhu.}

\FThm{\TaggThm{TianFullThm}
\ref{T3,T4,TZ}
Let $(M,\JJJ)$ be a Fano manifold and 
$G$ be a compact subgroup of $\Aut(M,\JJJ)$.
If $F$ or $E_0$ is proper on $\HcGplus$ then 
$(M,\JJJ)$ admits a $G$-invariant \KE metric.
Conversely, if $(M,\JJJ)$ admits a $G$-invariant \KE metric
then $F$ and $E_0$ are proper on $\Hconeplus$.
}

\noindent

We remark that when $\aut(M,\JJJ)$ is semisimple then $\HcG\sseq\Hc(\Lambda_1)$ \ref{PSSW}.

Let us turn to the proof of our main theorems and begin with \Ref{PropernessThm}.
Assume that a \KE form $\o$ exists. Then $F$ is proper on $\Hcplus$
by \Ref{TianThm}. By \Ref{DTLemma} and \fovpNegativeeq\ so is $E_0$.
From \Ref{BMRProp} we have
$$
E_{k+1}(\o,\ovp)
=  
E_{k}(\o,\ovp)
+(I_{k+1}-I_{k})(\ovp,\Ricovp)-(I_{k+1}-I_{k})(\o,\Rico),
$$
with
$I_{k+1}\ge I_{k}$ on $\Hc\times\Hc$ as noted after \Ikineq.
It follows that if $E_{k}$ is proper on $\Hcplus$ so is $E_{k+1}$. 
We conclude that $E_n$ is proper on $\Hcplus$. 

Assume that $E_n$ is proper on $\Hcplus$.
Then from \Ref{FzeroEn} and the
Calabi-Yau Theorem we see 
that $F$ is bounded from below on $\Hc$ and from \Ref{DTLemma} and
\fovpNegativeeq\
it follows that so is $E_0$.
Therefore from \Ref{AubDef}, given $\o\in\Hc$, the continuity path 
\MASecondpathseq\
extends 
for all $t<1$.

From the properness and exactness of $E_n$ there exists a function $\tau_\o$
as in \Ref{SectionEnergyFunctionals} satisfying 
$E_n(\o_{\vp_0},\ovpt)\ge\tau_\o(I(\o,\ovpt))-E_n(\o,\o_{\vp_0}).$ 
Hence it suffices now to show that $E_n(\o_{\vp_0},\ovpt)$ is uniformly
bounded from above for all $t> t_0$ with $t_0$ depending only on $(M,\JJJ,\o)$.
We will then have that $I(\o,\ovpt)$ is uniformly bounded 
independently of $t\in[0,1)$. This will entail  
a uniform bound on $\n{\vpt}{\Linf}$ \ref{A3, Proposition 7.35; T4, Lemma 6.19} and
hence a uniform bound
on $\n{\vpt}{\Ctwobefoot(M,g_\o)}$ for some $\be\in(0,1)$ \ref{A1,Y}. 
By the continuity method arguments therein one then
concludes that a unique smooth solution exists at $t=1$ that is a \K potential for 
a \KE form.

In fact we will find such a $t_0$ depending only on $n$ for each $E_k$. 
The computation 
that follows involves expressions similar to those 
that figure in 
the work of Song and Weinkove; using
\Ref{BMRProp} considerably simplifies our calculations compared to the ones 
there.

Fix $\tau\in[0,1]$. First, from \AubPatheq\ and the definition of $E_0$ we have
$$\eqalignno{
E_0(\o_{\vp_0},\o_{\vp_\tau})
= &\ 
\int_{[0,\tau]}\ddt E_0(\o_{\vp_0},\o_{\vp_t})dt
\cr
= &\
\V\int_{M\times [0,\tau]}
(1-t)n\dot\vpt\i\ddbar\vpt\w\ovpt^{n-1}\w dt
\cr
= & 
-\int_{[0,\tau]}(1-t)\ddt(I-J)(\o,\o_{\vp_t})dt
\cr
= &  
-(1-\tau)(I-J)(\o,\o_{\vp_\tau})
\cr 
& +(I-J)(\o,\o_{\vp_0})
-\int_{[0,\tau]}(I-J)(\o,\o_{\vp_t})dt.&\eqnno\Ezeroteq
}$$
From \Ref{BMRProp}, \IJineq\ and \IkJeq\ we therefore conclude that there
exists a constant $c_\o$ depending only on $(M,\JJJ,\o)$ for which
$$
(n+1)E_k(\o_{\vp_0},\ovpta)\le
-(1-\tau)I(\o,\o_{\vp_\tau})+nI(\ovpta,\Ric\ovpta)+c_\o.\eqn\Ektauineq
$$
From \AubPatheq\  
$$
\eqano{
I(\ovpta,\Ric\ovpta)
& =  
(1-\tau)^2\V\intm\i\del\vpta\w\dbar\vpta
\w\sum_{l=0}^{n-1}\ovpta^{n-l-1}\w(\tau\ovpta+(1-\tau)\o)^l
\cr
& \!\!\!\!\!\!\!\! =  (1-\tau)^2\V\intm\i\del\vpta\w\dbar\vpta
\w\sum_{l=0}^{n-1}\sum_{j=0}^l{l\choose j}\tau^{l-j}(1-\tau)^j\ovpta^{n-j-1}\w\o^j
\cr
& \!\!\!\!\!\!\!\! =  (1-\tau)^2\V\intm\i\del\vpta\w\dbar\vpta
\w
\sum_{j=0}^{n-1}(1-\tau)^j
\sum_{l=j}^{n-1}{l\choose j}\tau^{l-j}\ovpta^{n-j-1}\w\o^j.
}$$
Note that 
$$
(1-\tau)^j
\sum_{l=j}^{n-1}{l\choose j}\tau^{l-j}
\le
(1-\tau)^j(n-1){{n-1}\choose j}.\eqn\coeffeq
$$
We may choose $t_1\in[0,1)$ depending only on $n$ in such a way that for all $\tau\in[t_1,1]$
the expression on the right hand side of \coeffeq\ is smaller than $n$ for each
$j=0,\ldots,n-1$.
We conclude that
$$
I(\ovpta,\Ric\ovpta)\le n(1-\tau)^2I(\o,\ovpta), \qq\all\,\tau\in[t_1,1).\eqn\IIkteq
$$
Returning to \Ektauineq\ we then see that 
$E_k(\o_{\vp_0},\ovpta)\le c_\o/(n+1)$ whenever $\tau\in[\max\{t_1,1-\frac1{n^2}\},1)$.
This concludes the proof of \Ref{PropernessThm}.\done

As a corollary of the proof we record the following fact. 
\FCor{\TaggCor{CorrBound}
Let $(M,\h{\rm J})$ be a Fano manifold.
If one of the functionals
$F,E_0,\ldots,E_n$ is bounded from below on $\Hcplus$ so are the rest.}

Combining \Ref{CorrBound} with \Ref{BDTThm} concludes the proof of \Ref{BoundednessThm}.\done

We end this section with several remarks.

\FRemm{
\TaggRemm{RemmFull}%
Our methods imply that the refined version of 
\Ref{TianThm} (\Ref{TianFullThm}\unskip) also extends to each of the functionals $E_k$.
}

\FRemm{
\TaggRemm{RemmBound}
Note that one may
state \Ref{CorrBound} with $\Hcplus$ replaced by $\Hc$ for $F,E_0$ and $E_1$.
Indeed,  
recall that once $F$ is bounded from below
on $\Hcplus$ so are each of the $E_k$ while a lower bound for $E_n$ on $\Hcplus$
implies a lower bound for $F$ on $\Hc$ (by \Ref{FzeroEn}) which, in turn, implies the same for $E_0$
(using \Ref{DTLemma}) and for $E_1$ (using \Ref{BMRProp}).
Some special cases of \Ref{CorrBound} appeared previously, namely the fact
that when $F$ is bounded from below so is $E_0$ \ref{DT} and vice versa \ref{L2}, 
and the fact that when $E_0$ is bounded from below so is $E_1$ \ref{P} and vice versa \ref{CLW}.
}

\FRemm{
\TaggRemm{RemmLBounds}
Assume that the functionals $F$ and $E_k,\; k\in\{0,\ldots,n\}$ are bounded from below
on $\Hcplus$ and for each $\o\in\Hc$ set
$
l(\o)=\inf_{\ovp\in\Hcfoot} F(\o,\ovp)
$
and
$$
l_k(\o)=\cases{\inf_{\ovp\in\Hcfoot} E_k(\o,\ovp),& for  $k=0,1$,\cr
\inf_{\ovp\in\Hcfootplus} E_k(\o,\ovp),& for  $k=2,\ldots,n$.\cr
}
$$
Then the following relations hold between the various lower bounds: 
$$
l(\o)+\V\intM f_\o\o^n= l_0(\o)= l_k(\o)+I_k(\o,\Ric\o).\eqn
$$
This generalizes the relation between $l$ and $l_0$ \ref{L2} and between $l$ and $l_1$ \ref{CLW} 
that appeared recently; our proof, given below, appears considerably simpler.\hfill\break
\smallskip
\indent\proof
By \Ref{DTLemma}, \fovpNegativeeq\ and \Ref{BMRProp}
$$
l(\o)+\V\intM f_\o\o^n\le l_0(\o)\le l_k(\o)+I_k(\o,\Ric\o).\eqn\llkeq
$$
(For the second inequality we used \SecondBMRPropEq\ and the fact that $I_k(\ovp,\Ric\ovp)\ge0$ 
for $\ovp\in\Hcplus$.)
On the other hand, 
note first that from \Ezeroteq\ it follows that $\int_{[0,1]}(I-J)(\o,\ovpt)dt$ is bounded. As remarked 
in \Ref{SectionContinuity} the function $(I-J)(\o,\ovpt)$ is nondecreasing in $t$. 
Hence
$$
(1-\tau)(I-J)(\o,\ovpta)\le\int_{[\tau,1]}(I-J)(\o,\ovpt)dt,
$$
and therefore \ref{DT, p. 67}
$$
\lim_{\tau\ra1^-}(1-\tau)(I-J)(\o,\ovpta)=0.\eqn\IminusJteq
$$
Going back to \Ezeroteq\ and using the identity 
$E_0(\o,\o_{\vp_0})+(I-J)(\o,\o_{\vp_0})=\Vm\intm f_\o\o^n$ we have 
$$
\lim_{\tau\ra1^-}
E_0(\o,\ovpt)=\V\intm f_\o\o^n-\int_{[0,1]} (I-J)(\o,\ovpt)dt.
$$
By a theorem of Ding and Tian we have \ref{DT, Theorem 1.2}
$$
l(\o)=\lim_{t\ra1^-}F(\o,\ovpt)=-\int_{[0,1]} (I-J)(\o,\ovpt)dt.\eqn\lzeroeq
$$
Combining with \llkeq\ we conclude that 
$$
l_0(\o)=\lim_{t\ra1^-}E_0(\o,\ovpt)=l(\o)+\V\intM f_\o\o^n.
$$
Finally, using 
\IkJeq, \IJineq, \IIkteq\ and \IminusJteq\ it follows that
$\lim_{t\ra1^-}I_k(\ovpt,\Ric\ovpt)=0$. Therefore, using \Ref{BMRProp}
\SecondBMRPropEq\ again
we have $l_0(\o)\ge l_k(\o)+I_k(\o,\Ric\o)$. \done
}

\FRemm{
\TaggRemm{RemmProper}
Note that from \Ref{BMRProp} it follows that if $F$ is proper on $\Hc$ 
(equivalently on $\Hcplus$) with 
$F(\o,\ovp)\ge \tau_\o((I-J)(\o,\ovp))$ then we have the inequality
$E_k(\o,\ovp)\ge \tau_\o((I-J)(\o,\ovp))-I_k(\o,\Ric\o)$ on $\Hcplus$ 
(and for $k=0,1$ 
on $\Hc$). On the determination of explicit functions $\tau_\o$ we refer to \ref{PSSW,T3,T4}.
}

\subsection{Boundedness of energy functionals and the Moser-Trudinger-Onofri inequality}
\TaggSection{SectionMTO}
In this section we suppose that a \KE metric $\o$ exists.
First, we state the following fundamental theorem:

\FThm{\TaggThm{BMBound}\ref{BM, Theorem A, Corollary 8.3; Ba, Theorem 1} 
Let $(M,\JJJ,\o)$ be a \KE Fano manifold. Then
$E_0(\o,\ovp)\ge0$ for all $\ovp\in\Hc$ and $E_n(\o,\ovp)\ge0$ for all $\ovp\in\Hcplus$
with equality if and only
if $\ovp=h^\star\o$ with $h\in\Aut(M,\JJJ)_0$.
}

Building on these results, Song and Weinkove proved: (i) The first statement holds with $E_0$ replaced by $E_1$  
(see also \ref{P}), and (ii) the second statement holds with $E_n$ replaced by $E_k$ for each
$k\in\{2,\ldots,n-1\}$. \Ref{BMRProp} provides a much simplified proof of these two facts. Moreover, it allows
to improve on (ii). Let 
$$
\calA_k:=\{\ovp\in\Hc\,:\, E_k(\o,\ovp)\ge0\}.\eqn
$$
Then we have shown that
$$
\calA_k\supseteq\calB_k:=\{\ovp\in\Hc\,:\, I_k(\ovp,\Ricovp)\ge0\}.\eqn
$$
For example, for $k=1$ this gives $\calA_1=\Hc$, when $k=2$ we have
$$
\calA_2\supseteq\calB_2\supseteq\{\ovp\in\Hc\,:\, \Ricovp+2\ovp\ge0\},
$$
for $k=3$ 
$$
\calA_3\supseteq\calB_3\supseteq\{\ovp\in\Hc\,:\, \Ricovp+\ovp\ge0\},
$$
and for arbitrary $k$ one may readily obtain an explicit bound (depending on $k$) on the set
$\calB_k$, and hence on 
$\calA_k$, in terms
of a lower bound on the Ricci curvature, using the definition \Ikeq.

Let $\oFSc$ denote the Fubini-Study form of constant Ricci curvature $c$ on $(S^2,\JJJ)$, the Riemann sphere, given locally by
$$
\oFSc=\frac\i{c\pi}\frac{dz\w\dbz}{(1+|z|^2)^2}.
$$
Here $V=\int_{S^2}\oFSc=c_1([M])/c=2/c$. For
$c=1/2\pi$ it is induced from restricting the Euclidean metric on $\RR^3$ to 
the radius $1$ sphere. Denote by $W^{1,2}(S^2)$ the space of functions on $S^2$
that are square-summable and so is their gradient (with respect to some Riemannian 
metric).
The Moser-Trudinger-Onofri inequality states:
\FThm{\ref{Mo,O,Tr}
\TaggThm{MTOThm}
 For
$\o=
\o_{\hbox{\small FS},2/V}
$ and 
any function $\vp$ on $S^2$ in $W^{1,2}(S^2)$ one has
$$
\V\int_{S^2} e^{-\vp+\V\int_{S^2}\vp\o}\o\le e^{\V\int_{S^2}\ha\i\del\vp\w\dbar\vp}.\eqn\MTClassicalineq
$$
Equality holds if and only if $\ovp$ is the pull-back of $\o$ by a M\"obius transformation.}
\noindent
An alternative proof of this inequality has been given by Ding and Tian for those functions $\vp$ that
belong to the subspace 
$\calH_\o\sseq W^{1,2}(S^2)$. The proof uses the properties of $F$.
We now note that our work provides a new and succinct proof of the original
Moser-Trudinger-Onofri inequality entirely within the framework of exact energy functionals.
This is the first proof that does not use symmetrization/rearrangement arguments. Other proofs
of this inequality have been given by Onofri \ref{O}, Hong \ref{H}, Osgood-Phillips-Sarnak \ref{OPS},
Beckner \ref{B}, Carlen and Loss \ref{CL1,CL2}, Ghigi \ref{G} (for more background we refer to Chang \ref{Ch}).

\medskip
\proof
By \Ref{BMBound} and \Ref{BMRProp} $E_1(\o,\cdot)\ge0$ on $\Ho$. 
Given $\vp\in\Cinf(S^2)$ there exists $\psi\in\calH_\o$ such that
$\Ric\o_\psi=\ovp$ by solving the Poisson equation on $S^2$.
Thus by \Ref{FzeroEn} $F(\o,\cdot)\ge0$ on $\Do$. Using the definition of $F$, 
for any smooth function $\vp$ we obtain \MTClassicalineq. Since 
$\Cinf(S^2)$ is dense in $W^{1,2}(S^2)$ we conclude.\done

Ding and Tian showed that a restricted analogue of this inequality holds also for higher dimensional %
manifolds:
\FThm{
\TaggThm{DTBoundednessThm}
\ref{DT}
Let $(M,\JJJ)$ be a Fano manifold and let $\o\in\Hc$.
Assume that $F$ is bounded from below on $\Hc$
and let $a=-\inf_{\hbox{\calfoot H}_{c_1}} F(\o,\cdot)$. Then 
for each $\vp\in\H_\o$ holds
$$
\V\intm e^{-\vp+\V\intm\vp\on}\on
\le 
e^{J(\o,\ovp)+a}.\eqn\GeneralizedMTOeq
$$
If $(M,\JJJ,\o)$ is \KE then $a=0$.
}

Recall Jensen's inequality $\V\intm e^{-\vp+\V\intm\vp\on}\on\ge1$ \ref{HLP}.
Now observe that in higher dimensions, due to \Ref{UnboundedLemma}, inequality \GeneralizedMTOeq\ cannot be
extended to all of $\CinfM$. A natural question is therefore: On a \KE manifold, 
what is the largest neighborhood of $\H_\o$ inside $\CinfM$ on which \GeneralizedMTOeq\ does hold?
Naturally, we call such a neighborhood the
Moser-Trudinger-Onofri neighborhood of $\H_\o$, and put
$$
MTO_n=\{\vp\in\CinfM: \!\vp \hbox{$\!$ satisfies \GeneralizedMTOeq\ on the Fano manifold\ } 
(M,\JJJ), \dim_\CCfoot M=n\}.
\eqn\MTOnEq
$$ 

Using \Ref{FzeroEn}, we have the following characterization of the Moser-Trudinger-Onofri neighborhood.
By abuse of notation we do not distinguish here between the set $\Ric(\calA_n)$ in $\H_\O$ and the corresponding
set in $\H_\o$.

\FThm{
\TaggThm{CharacterizationThm}
Let $(M,\JJJ,\o)$ be a \KE Fano manifold. Then $\vp\in\CinfM$ satisfies the generalized
Moser-Trudinger-Onofri inequality \GeneralizedMTOeq\ if and only if
there exists a function $\psi\in\CinfM$ such that $\Ric\o_\psi=\ovp$
and $\o_\psi\in\calA_n$. That is, $MTO_n=\Ric(\calA_n)\supset\H_\o.$
}

Recall that $\calB_n\subseteq\calA_n$ and that we have bounds on $\calB_n$ in terms of the Ricci
curvature. Therefore, \Ref{CharacterizationThm} shows that in higher dimensions the Moser-Trudinger-Onofri
inequality is related to Ricci curvature and holds on a set strictly larger than the space of \K potentials.
It would be interesting to improve the bounds both on $\calA_n$ and on $\calB_n$.

We now state another corollary of our arguments.

\FCor{
\TaggCor{EnunboundedCor} Let $(M,\JJJ,\o)$ be a Fano manifold. Then the
Ricci energy $E_n$ is unbounded from below on $\Hc$ if and only if $n>1$.}

Before concluding, we remark that 
\Ref{DTBoundednessThm} can be strengthened
using the results obtained here. 
The same applies to later extensions of this inequality \ref{PSSW,T4} and will 
figure in a subsequent article.

\subsectionno{Appendix}
In this appendix we prove \Ref{BMRProp}.
First, in order to establish formula \SecondBMRPropEq\ we show that the variations of both sides of the equation agree. 

$$\leqalignno{
-(k+1)V\ddt I_k(\ovp,\Ricovp)
& =
\ddt \intm \fovp\i\ddbar\fovp\w\sum_{l=0}^{k-1} (k-l)\ovp^{n-1-l}\w(\Ricovp)^l
\cr
& =
\ddt \intm\! \fovp(\Ricovp-\ovp)\!\w\!\sum_{l=0}^{k-1} (k-l)\ovp^{n-1-l}\w(\Ricovp)^l
\cr
& =
\ddt \intm \fovp\Big(-k\ovpn+\sum_{l=1}^{k}\ovp^{n-l}\w(\Ricovp)^l\Big)
\cr
& = \intm\dot\fovp
	\Big(-k\ovpn+\sum_{l=1}^{k}\ovp^{n-l}\w(\Ricovp)^l\Big)&\eqnno\ITermeq
\cr
  + \intm &\fovp\i\ddbar\dot\vp\w
        \Big(\sum_{l=1}^{k}(n-l)\ovp^{n-l-1}\w(\Ricovp)^l-kn\ovp^{n-1}\Big)&\eqnno\IITermeq
\cr
& - \intm\fovp\i\ddbar\Dvp\dot\vp
	\sum_{l=1}^{k} l\ovp^{n-l}\w(\Ricovp)^{l-1}.&\eqnno\IIITermeq
}$$
First, we write \ITermeq\ as
$$
\intm\dot\fovp
     \Big(-k\ovpn+\sum_{l=1}^{k}\ovp^{n-l}\w(\Ricovp)^l\Big)
=:\iota_1+\mu_1.
$$
We will evaluate \IITermeq\ and \IIITermeq\ by substituting once again 
$\i\ddbar\fovp=\Ricovp-\ovp$.
For \IITermeq\ we get 
$$
\eqano{
&
\intm\dot\vp(\Ricovp-\ovp)\w
        \Big(-kn\ovp^{n-1}+\sum_{l=1}^{k}(n-l)\ovp^{n-l-1}\w(\Ricovp)^l\Big)
\cr  
& \q =
\intm\dot\vp
        \Big(kn\ovp^{n}-kn\ovp^{n-1}\w\Ricovp
\cr
&\qq\qq\qq -(n-1)\ovp^{n-1}\w\Ricovp+\sum_{l=2}^{k}\ovp^{n-l}\w(\Ricovp)^{l}
\cr
&\qq\qq\qq     +(n-k)\ovp^{n-k-1}\w(\Ricovp)^{k+1}\Big)
\cr  
& \q =
\intm\dot\vp
        \Big([-(n-k)+(k+1)n-k]\ovp^{n}-(k+1)n\ovp^{n-1}\w\Ricovp
\cr
&\qq\qq\qq
+\sum_{l=1}^{k}\ovp^{n-l}\w(\Ricovp)^{l}
     +(n-k)\ovp^{n-k-1}\w(\Ricovp)^{k+1}\Big)
\cr
&\q =:(\kappa_1+\lambda_1+\iota_2)+\lambda_2+\mu_2+\kappa_2.
}$$

For \IIITermeq\ we get
$$
\eqano{
&\intm\Dvp\dot\vp(\ovp-\Ricovp)\w
	\sum_{l=1}^{k} l\ovp^{n-l}\w(\Ricovp)^{l-1}
\cr 
& \q =
\intm\Dvp\dot\vp\Big(\ovpn+
	\sum_{l=1}^{k-1} \ovp^{n-l}\w(\Ricovp)^{l}
        -k\ovp^{n-k}\w(\Ricovp)^k
\Big)
\cr 
& \q =
\intm\Dvp\dot\vp\Big(
	\sum_{l=1}^{k} \ovp^{n-l}\w(\Ricovp)^{l}
        -(k+1)\ovp^{n-k}\w(\Ricovp)^k
\Big)
\cr 
& \q =:\mu_3+\kappa_3.
}$$

Noting that $\dot\fovp=-\Dvp\dot\vp-\dot\vp+c$ with $c$ a constant
yields $\iota_1+\iota_2=-kcV$ and $\mu_1+\mu_2+\mu_3=kcV$.
Note that $\kappa_1+\kappa_2+\kappa_3=-(k+1)V\ddt E_k(\o,\ovp)$
and $\lambda_1+\lambda_2=(k+1)V\ddt E_0(\o,\ovp)$. This completes the proof
of \SecondBMRPropEq. 

Formulas \FirstBMRPropEq\ and \ThirdBMRPropEq\ 
now follow: first use \SecondBMRPropEq\ with $k=n$ to express
$E_0$ in terms of $E_n$ and $J$, and then substitute this expression back
into \SecondBMRPropEq\ and apply \Ref{IkJkLemma}. \done

Let us note that one way one could arrive at the formula
would be to use the expression for $(k+1)I_k-kI_{k-1}$ (see \Ikeq) 
and  \Ref{DTLemma} together with the observation
$$
\ddt ((k+1)E_k-kE_{k-1})(\o,\ovpt)=-\V\intm\dot\vpt\ovptn
-
\ddt\Big(\V\intM f_{\ovpt}(\Ricovpt)^k\w\ovpt^{n-k}\Big).\eqn\Ekdifferenceeq
$$

\def\smallblackbox{\vrule height.6ex width .5ex depth -.1ex}
\def\boxseparation{\hfil\smallblackbox$\q$\smallblackbox$\q$\smallblackbox\hfil}
\bigskip
\boxseparation
\bigskip

I would like to express my deep gratitude to my teacher, Gang Tian.
I thank G. Maschler for his interest in this work, N. Pali for a useful discussion, 
J. Song and V. Tosatti for helpful discussions as well as useful comments,
and a referee for a careful reading of this manuscript.
I thank William Browder for his kindness and my 908 Fine Hall office mates for their pleasant
company. This material is based upon work supported under a National Science 
Foundation Graduate Research Fellowship.

\frenchspacing

\bigskip\bigskip
\noindent{\bf Bibliography}
\bigskip
\def\ref#1{\Taggf{#1}\item{ {\bf[}{\sans #1}{\bf]} } }

\ref{Al} Lars Alexandersson, On vanishing-curvature extensions of Lorentzian metrics,
{\sl The Journal of Geometric Analysis} {\bf 4} (1994), 425--466.
\sm
\ref{A1} Thierry Aubin, \'{E}quations du type {M}onge-{A}mp\`ere sur les vari\'et\'es
              k\"ahl\'eriennes compactes, {\sl Bulletin des Sciences Math\'ematiques} 
{\bf 102} (1978), 63--95.
\sm
\ref{A2} \opcit, R\'eduction du cas positif de l'\'equation de
              {M}onge-{A}mp\`ere sur les vari\'et\'es k\"ahl\'eriennes
              compactes \`a la d\'emonstration d'une in\'egalit\'e, 
{\sl Journal of Functional Analysis} {\bf 57} (1984), 143--153.
\sm
\ref{A3} \opcit, Some Nonlinear Problems in Riemannian Geometry, Springer, 1998.
\sm
\ref{Ba} Shigetoshi Bando, The K-Energy Map, almost \KE metrics and an
inequality of the Miyaoka-Yau type, {\sl T\^ohoku Mathematical Journal}
{\bf 39} (1987), 231--235.
\sm
\ref{BM} Shigetoshi Bando, Toshiki Mabuchi, Uniqueness of \KE metrics
modulo connected group actions, in {\it Algebraic Geometry,
Sendai, 1985} (T. Oda, Ed.), Advanced Studies in Pure Mathematics {\bf 10},
Kinokuniya, 1987, 11--40.
\sm
\ref{B} William Beckner, Sharp Sobolev inequalities on the sphere and the 
Moser-Trudinger inequality, {\sl Annals of Mathematics} {\bf 138} (1993), 213--242.
\sm
\ref{Be} Arthur L. Besse, Einstein manifolds, Springer, 1987.
\sm
\ref{Bi} Olivier Biquard, M\'etriques k\"ahl\'eriennes \`a courbure 
scalaire constante: Unicit\'e, stabilit\'e, preprint. 
To appear in {\sl Ast\'erique}.
\sm
\ref{C} Eugenio Calabi, Extremal {K}\"ahler metrics, {II}, in
{\it Differential geometry and complex analysis} (I. Chavel, H. M. Farkas, Eds.), Springer, 1985, {95--114}.
\sm
\ref{CL1} Eric A. Carlen, Michael Loss, Competing symmetries of some functionals arising in mathematical physics, 
in {\it Stochastic processes, physics and geometry} (S. Albeverio et al., Eds.), World
Scientific, 1990, 277--288.
\sm
\ref{CL2} \opcit, Competing symmetries, the logarithmic HLS inequality and Onofri's
inequality on $S^n$, {\sl Geometric and Functional Analysis} {\bf 2} (1992), 90--104.
\sm
\ref{Ch} Sun-Yung A. Chang, Non-linear elliptic equations in conformal geometry,
European Mathematical Society, 2004.
\sm
\ref{CLW} Xiu-Xiong Chen, Hao-Zhao Li, Bing Wang,
On the K\"ahler-Ricci flow with small initial $E_1$ energy (I),
 preprint, arxiv: math.DG/0609694 v2. To appear in {\sl Geometric and Functional Analysis}.
\sm
\ref{CT} Xiu-Xiong Chen, Gang Tian, Ricci flow on \KE surfaces,
{\sl Inventiones Mathematicae} {\bf 147} (2002), 487--544.
\sm
\ref{D} Wei-Yue Ding, Remarks on the existence problem of positive
              {K}\"ahler-{E}instein metrics, {\sl Mathematische Annalen} {\bf 282} (1988), 463--471.
\sm
\ref{DT} Wei-Yue Ding, Gang Tian, The generalized Moser-Trudinger inequality, 
in {\it Nonlinear Analysis and Microlocal Analysis: Proceedings of the International
Conference at Nankai Institute of Mathematics} (K.-C. Chang et al., Eds.), 
World Scientific, 1992, 57--70. ISBN 9810209134.
\sm
\ref{Do} Simon K. Donaldson, Anti self-dual Yang-Mills connections over complex
                     algebraic surfaces and stable vector bundles, {\sl Proceedings
                     of the London Mathematical Society} {\bf 50} (1985), 1--26.
\sm
\ref{Fu} Akira Fujiki, On automorphism groups of compact \K manifolds,
{\sl Inventiones Mathematicae} {\bf 44} (1978), 225--258.
\sm
\ref{F1} Akito Futaki, K\"ahler-{E}instein metrics and integral invariants,
Lecture Notes in Mathematics {\bf 1314}, Springer, 1988.
\sm
\ref{F2} \opcit, Stability, integral invariants and canonical \K metrics, preprint, 2005.
\sm
\ref{F3} \opcit, Asymptotic Chow semi-stability and integral invariants, 
{\sl International Journal of Mathematics} {\bf 15} (2004), 967--979. 
\sm
\ref{G} Alessandro Ghigi, On the Moser-Onofri and Pr\'ekopa-Leindler inequalities,
{\sl Collectanea Mathematica} {\bf 56} (2005), 143--156.
\sm
\ref{HLP} Godfrey H. Hardy, John E. Littlewood, George P{\'o}lya,
           Inequalities (Second Edition), Cambridge University Press, 1952.
\sm
\ref{H} Chong-Wei Hong, A best constant and the Gaussian curvature, {\sl Proceedings of the 
American Mathematical Society} {\bf 97} (1986), 737--747.
\sm
\ref{L1} Hao-Zhao Li, A new formula for the Chen-Tian energy functionals $E_k$ and
its applications, preprint, arxiv: math.DG/0609724 v1.
\sm
\ref{L2} \opcit, On the lower bound of the K-energy and $F$ functional, preprint,
arxiv: math. DG/0609725 v1.
\sm
\ref{Li} Chiung-Ju Liu, Bando-Futaki invariants on hypersurfaces,
preprint, arxiv: math. DG/0406029 v3.
\sm
\ref{M} Toshiki Mabuchi, K-energy maps integrating Futaki invariants,
{\sl T\^ohoku Mathematical Journal} {\bf 38} (1986), 575--593.
\sm
\ref{Ma} Gideon Maschler, Central \K metrics, 
{\sl Transactions of the American Mathematical Society} {\bf 355}
(2003), 2161--2182.
\sm
\ref{Mat} Yoz{\^o} Matsushima, Sur la structure du groupe d'hom\'eomorphismes analytiques
              d'une certaine vari\'et\'e k\"ahl\'erienne, {\sl Nagoya Mathematical Journal} {\bf 11}
              (1957), 145--150.
\sm
\ref{Mo} J\"urgen Moser, A sharp form of an inequality by {N}. {T}rudinger,
{\sl Indiana University Mathematics Journal} {\bf 20} (1971), 1077--1092.
\sm
\ref{O} Enrico Onofri, On the positivity of the effective action in a theory of
random surfaces, {\sl Communications in Mathematical Physics} {\bf 86} (1982), 321--326.
\sm
\ref{OPS} Brad Osgood, Ralph Phillips, Peter Sarnak, Extremals of determinants of {L}aplacians,
{\sl Journal of Functional Analysis} {\bf 80} (1988), 148--211.
\sm
\ref{P} Nefton Pali, A consequence of a lower bound of the K-energy,
  {\sl International Mathematics Research Notices} (2005), 3081--3090.
\sm
\ref{PSSW} 
Duong H. Phong, Jian Song, Jacob Sturm, Ben Weinkove,
The Moser-Trudinger inequality on \KE manifolds, preprint,
arxiv: math.DG/0604076 v2.
\sm
\ref{R} Yanir A. Rubinstein, Ph.D. thesis, Massachusetts Institute of Technology. In preparation.
\sm
\ref{S1} Yum-Tong Siu,
The existence of {K}\"ahler-{E}instein metrics on manifolds
with positive anticanonical line bundle and a suitable finite symmetry group,
{\sl Annals of Mathematics} {\bf 127} (1988),  585--627.
\sm
\ref{S2} \opcit, Lectures on Hermitian-Einstein metrics for stable bundles
and \KE metrics,  Birkh\"auser, 1987. 
\sm
\ref{SW} Jian Song, Ben Weinkove, Energy functionals and canonical \K metrics,
{\sl Duke Mathematical Journal} {\bf 137} (2007), 159--184.
\sm
\ref{Sz} G\'abor Sz\'ekelyhidi, Extremal metrics and K-stability, Ph.D. Thesis, Imperial College, 
2006. Available at arxiv: math.DG/0611002 v1.
\sm
\ref{Th} Richard P. Thomas, 
Notes on GIT and symplectic reduction for bundles and varieties, 
in {\it Surveys in Differential Geometry: Essays in memory of S.-S. Chern}
(S.-T. Yau, Ed.), International Press, 2006, 221--273.
\sm
\ref{T1} Gang Tian, 
On K\"ahler-Einstein metrics on certain K\"ahler manifolds with\break $c_1(M)>0$,
{\sl Inventiones Mathematicae} {\bf 89} (1987), 225--246.
\sm
\ref{T2} \opcit,
On Calabi's conjecture for complex surfaces with positive first Chern class,
{\sl Inventiones Mathematicae} {\bf 101} (1990), {101--172}.
\sm
\ref{T3} \opcit,
K\"ahler-{E}instein metrics with positive scalar curvature, 
{\sl Inventiones Mathematicae} {\bf 130} (1997), {1--37}.
\sm
\ref{T4} \opcit,
Canonical Metrics in \K Geometry, Birkh\"auser, 2000.
\sm
\ref{TY} Gang Tian, Shing-Tung Yau,
K\"ahler-Einstein metrics on complex surfaces with $c_1>0$,
{\sl Communications in Mathematical Physics} {\bf 112} (1987), 175--203.
\sm
\ref{TZ} Gang Tian, Xiao-Hua Zhu, A nonlinear inequality of Moser-Trudinger type,
{\sl Calculus of Variations} {\bf 10} (2000), 349--354.
\sm
\ref{To} Valentino Tosatti, On the critical points of the $E_k$ functionals
in \K geometry, preprint, arxiv: DG/0506021 v1. To appear in {\sl Proceedings of the 
American Mathematical Society}.
\sm
\ref{Tr} Neil S. Trudinger, On imbeddings into {O}rlicz spaces and some applications,
{\sl Journal of Mathematics and Mechanics} {\bf 17} (1967), 473--483.
\sm
\ref{Y} Shing-Tung Yau, On the Ricci curvature of a compact \K
manifold and the Complex \MA equation, I, {\sl Communications in Pure
and Applied Mathematics} {\bf 31} (1978), 339--411.
\end

\bigskip
\bigskip
\noindent
{\caps Yanir A. Rubinstein}
\smallskip
\medskip

\noindent
{\caps Address: 908 Fine Hall, Princeton University, Princeton, NJ 08544, USA}

\noindent
{\caps Email:} {\tt yanir@math.princeton.edu}, {\tt yanir@member.ams.org}

\end